\documentclass[10pt]{amsart}

\usepackage[all, cmtip]{xy}
\usepackage{amsmath,amsthm,amssymb,amsfonts,amscd}
\usepackage{tikz-cd}
\usepackage{float} 

\usepackage{graphicx} 
 
\usepackage{kotex}

\usepackage{amsmath,amscd}
\usepackage[all,cmtip]{xy}

\newtheorem{theorem}{Theorem}[section]
\newtheorem{lemma}[theorem]{Lemma}
\newtheorem{proposition}[theorem]{Proposition}
\newtheorem{corollary}[theorem]{Corollary}

\theoremstyle{definition}
\newtheorem{definition}[theorem]{Definition}

\newtheorem{conjecture}[theorem]{Conjecture}

\theoremstyle{remark}
\newtheorem{remark}[theorem]{Remark}

\numberwithin{equation}{section}

\begin{document}

\title[Algebraic Montgomery-Yang problem and cascade conjecture]{Algebraic Montgomery-Yang problem \\ and cascade conjecture}

\author{DongSeon Hwang}
\address{Department of Mathematics, Ajou University, Suwon 16499, Republic Of Korea}
\email{dshwang@ajou.ac.kr}


\subjclass[2010]{Primary 14J25; Secondary 14J17, 14J26.}
 
\date{January 10, 2021}
 
\keywords{$\mathbb{Q}$-homology projective plane, algebraic Montgomery-Yang problem, cascade,  $\mathbb{P}^1$-fibration} 

\begin{abstract}
The conjecture called algebraic Montgomery-Yang problem is still open for rational $\mathbb{Q}$-homology projective planes with cyclic  quotient singularities having ample canonical divisor. All known such surfaces have a special birational behavior 
called a cascade. In this note, we establish algebraic Montgomery-Yang problem assuming  the cascade conjecture, which claims that every  rational $\mathbb{Q}$-homology projective planes with quotient singularities having ample canonical divisor admits a cascade. 
\end{abstract}

\maketitle

\section {Introduction}

Motivated by works of Seifert(\cite{Sei}), Montgomery and Yang(\cite{MY}) and Petri(\cite{Pet}), Fintushel and Stern formulated the following conjecture. 
\begin{conjecture}\cite{FS87} (Montgomery-Yang problem) Every pseudo-free differentiable circle action on the $5$-dimensional sphere has at most three non-free orbits.
\end{conjecture}

One way to approach this conjecture, as was already considered in \cite{MY}, \cite{FS85} and \cite{FS87}, is to consider the orbit space as a $4$-dimensional orbifold. Considering the  case when the orbit space has a structure of a complex projective surface,  Koll\'ar formulated the following version of the conjecture.

\begin{conjecture}\cite{Kol08}\label{amy} (Algebraic Montgomery-Yang problem) \\
Let $S$ be a $\mathbb{Q}$-homology projective plane with quotient singularities. If the smooth locus of $S$ is simply-connected, then $S$ has at most three singular points.
\end{conjecture}

See Notation \ref{notation} for the definition of a $\mathbb{Q}$-homology projective plane. 
The conjecture  turned out to be true when $S$ has a non-cyclic singular point(\cite{HK2011b}), $S$ is not a rational surface(\cite{HK2013}), and $-K_S$ is nef(\cite{HK2014} and \cite[Lemma 3.6 (4)]{HK2013}). Thus, it remained open only for rational $\mathbb{Q}$-homology projective planes with at worst   cyclic quotient singularities having  ample canonical divisor.

Several efforts have been devoted to construct rational $\mathbb{Q}$-homology projective planes  with ample canonical divisor (\cite{KM}, \cite{Kol08}, \cite{HK2012}, \cite{AL1}, \cite{AL2}). In fact, all of them admit  a  special  birational behavior, called a cascade.

\begin{definition}
Let $S$ be a $\mathbb{Q}$-homology projective plane with quotient singularities. We say that  $S$ \emph{admits a cascade}  if there exists a diagram as follows:

 $$\begin{CD}
S' = S'_t  @>\phi_t >> S'_{t-1}  @>\phi_{t-1} >> \ldots      @>\phi_1>> S'_{0}    \\ 
@V\pi_tVV @V\pi_{t-1}VV @.     @V\pi_0VV \\
S_t   := S @. S_{t-1} @. \ldots     @. S_{0}
\end{CD}$$
 
where for each $k$ 
\begin{enumerate}
\item $\phi_{k}$ is a  blow-down,  
\item $\pi_k$ is the minimal resolution, 
\item  $S_k$ is a $\mathbb{Q}$-homology projective plane,  and 
\item $S_0$ is a  log del Pezzo surface of Picard number one.
\end{enumerate}  
In this case, we also say that $S$ \emph{admits a cascade to $S_0$}.  
\end{definition}

Intuitively speaking, the idea behind this definition is that  $K_{S_t}$ becomes more negative for each step of cascade, i.e., as $t$ decreases. Thus, the study on rational $\mathbb{Q}$-homology projective plane with ample canonical divisor can be reduced to that on log del Pezzo surfaces of Picard number one.

We would like to pose the following conjecture, which is supported by every known example (\cite{KM},\cite{HK2012}, \cite{AL1}, \cite{AL2}). 

\begin{conjecture}\label{cascadeconj}
Every rational $\mathbb{Q}$-homology projective plane with quotient singularities such that the canonical divisor is ample admits a cascade. 
\end{conjecture}

Conjecture \ref{cascadeconj} implies that  the rational $\mathbb{Q}$-homology projective plane with quotient singularities having ample canonical divisor constructed in \cite[Theorem 8.2]{AL1} attains the minimal volume, as was expected.

In addition, the main result of the paper is the following.

\begin{theorem}\label{cascade=>amy}
Conjecture \ref{cascadeconj} implies Conjecture \ref{amy}. 
\end{theorem}

In fact, we prove the following stronger statement. Note that if Conjecture \ref{cascadeconj} is true, then there exists a  $(-1)$-curve $E$ with  $E.\mathcal{D} \leq 2$ where $\mathcal{D}$ is the reduced exceptional divisor of the minimal resolution of $S$.

\begin{theorem}\label{maintheorem}
Under the assumption in Conjecture \ref{amy}, we further assume that the canonical divisor is ample and  there exists a  $(-1)$-curve $E$ with  $E.\mathcal{D} \leq 2$ where $\mathcal{D}$ is the reduced exceptional divisor of the minimal resolution of $S$. Then, $S$ has at most three singular points.
\end{theorem}

The proof consists of arguments used in  \cite{HK2013}  together with   a method of using $\mathbb{P}^1$-fibration structure.

We work over the field $\mathbb{C}$ of complex numbers.

\section{Preliminaries}

\subsection{Notation}\label{notation}
We first fix some notations, following \cite{HK2013}, that will be used in the remaining of the paper. A normal projective surface with quotient singularities is called a \emph{$\mathbb{Q}$-homology projective plane} if it has the same Betti numbers as the complex projective plane.  It is said to be \emph{rational} if it is birationally equivalent to the projective plane $\mathbb{P}^2$. It is also said to be \emph{of log general type} if its canonical divisor is ample. We always denote by  $S$  a rational $\mathbb{Q}$-homology projective plane of log general type with $4$ cyclic quotient singularities $p_1$, \ldots $p_4$, unless otherwise stated. Let $f: S' \rightarrow S$ be a minimal resolution of $S$. Denote by $\mathcal{D}_i$ the reduced part of the $f$-exceptional divisor  $f^{-1}(p_i)$ of $p_i$ and by $l_i$ the number of irreducible components of $\mathcal{D}_i$. Let  $\mathcal{D} = \mathcal{D}_1 + \ldots + \mathcal{D}_4$ and $L = l_1 + \ldots + l_4$.

Let   $D$  be a $\mathbb{Q}$-divisor whose support consists of a chain of  rational curves with the dual graph   of the form $\overset{-n_1}{\circ}-\ldots-\overset{-n_l}{\circ}$ where each weight $n_i$ denotes the self-intersection number of the corresponding irreducible component $D_i$. Then, we introduce the following notation. (See \cite[Section 2]{HK2013} for more extensive explanation about the notation.)
\begin{enumerate}
    \item We denote by $q_D$ the order of the local fundamental group of the singular point  obtained by contracting the divisor $D$.
    \item We denote by $u_{k,D}$ the order of the local fundamental group of the singular point  obtained by contracting the  first $k-1$ irreducible components of $D$.
    \item We denote by $v_{k,D}$ the order of the local fundamental group of the singular point  obtained by contracting the  last $l-k$ irreducible components of $D$.
    \item For convenience, we also denote $v_{1,D}$ by $q_{1, D}$    and $u_{l,D}$ by $q_{l,D}$. 
    \item $tr_D = -\overset{l}{\underset{i=1}{\sum}} D^2_i$.
\end{enumerate} 
We will  omit the subscript $D$ in (5) whenever there is no confusion in the context.

\subsection{Known results and easy consequences}
We first summarize known results mainly from \cite{HK2013}. 

\begin{theorem}\cite{HK2011a}\label{maximum}
Let $S$ be a $\mathbb{Q}$-homology projective plane with quotient singularities. Then, $S$ has at most $5$ singular points and it has exactly $5$ singular points if and only if the singularities  of $S$ are of type  $3A_1+2A_3$ and the minimal resolution of $S$ is an Enriques surface. In particular, if $S$ is rational, then $S$ has at most $4$ singular points. 
\end{theorem} 

\begin{corollary}\label{ED>=2}
Assume the situation in Notation \ref{notation}. If $E$ is a $(-1)$-curve on $S'$, then $E.\mathcal{D} \geq 2$.  Moreover, if $E.\mathcal{D} = 2$, then $E$ intersects $\mathcal{D}$ at two different points. 
\end{corollary}

\begin{proof}
Since $S$ has Picard number one, $E.\mathcal{D} \geq 1$. If $E.\mathcal{D} = 1$, then  by blowing up the intersection point of $E$ and $\mathcal{D}$, we get a minimal resolution of another rational $\mathbb{Q}$-homology projective plane with $5$ singular points, a contradiction by Theorem \ref{maximum}.

Consider the case  $E.\mathcal{D} = 2$. Assume that  $E$ intersects $\mathcal{D}$ at one point $p$ with multiplicity $2$. If $p$ is an intersection of two different irreducible componnents of $\mathcal{D}$, by blowing up $p$, we get a minimal resolution of another rational $\mathbb{Q}$-homology projective plane with $5$ singular points, a contradiction by Theorem \ref{maximum}. Now, we may assume that $E$ intersects an irreducible component $D$ of $\mathcal{D}$ with multiplicity two.  By blowing up the intersection point of $E$ and $D$, and then again by blowing up the intersection point of the proper transform of $E$ and $D$, we get a minimal resolution of another $\mathbb{Q}$-homology projective plane which has $6$ quotient singularities, a contradiction by Theorem \ref{maximum}.
\end{proof}

Now, we consider possible orders for the local fundamental groups.

\begin{lemma} \cite[Lemma 5.2 and Lemma 5.3]{HK2013}\label{first4-tuple}
Let $S$ be as in Notation \ref{notation}. Assume that $H_1(S^{sm}, \mathbb{Z})$ is trivial. Then, the orders of the local fundamental groups of each singular point is either $(2,3,7,19)$ or $(2,3,5,q)$ for some positive integer $q$ with $gcd(q, 30)=1$. Moreover, the singularity type must be $A_1+A_2+\frac{1}{7}(1,1)+\frac{1}{19}(1,9)$ in the first case and the order $3$ singularity must be of type $\frac{1}{3}(1,1)$  in the latter case.
\end{lemma}

We have more precise information on the configuration of $(-1)$-curves on $S'$ in the first case.

\begin{lemma}\cite[Lemma 5.6]{HK2013}\label{23719conf}
Let $S$ be a rational $\mathbb{Q}$-homology projective plane of log general type  with exactly four cyclic quotient singular points $p_1, p_2, p_3, p_4$ of orders $(2, 3, 7, 19)$.   Let $\mathcal{D}$ be the reduced $f$-exceptional divisor on $S'$, and $E$ be any $(-1)$-curve on $S'$.  Then,  $E.\mathcal{D} \geq 2$ and the equality holds if and only if  $E.f^{-1}(p_i) = 0$ for $i = 1, 2, 3$, $E.f^{-1}(p_4) = 2$ and $E$ does not meet an end component of $f^{-1}(p_4)$.  
\end{lemma}

We also have more   information in the latter case.

\begin{lemma}\cite[Lemma 5.4 and Lemma 5.5]{HK2013}\label{235type}
Let $S$ be a rational $\mathbb{Q}$-homology projective plane of log general type  with exactly four cyclic quotient singular points $p_1, p_2, p_3, p_4$ of orders $(2, 3, 5, q)$.  Assume that the order $3$ singularity is
of type $\frac{1}{3}(1,1)$. Then, we have
\begin{enumerate}
    \item $L \geq 11$ and   $L = 11$ if and only if $S$ has singularities of one of the following types:
    \begin{enumerate}
        \item $A_1 + \frac{1}{3}(1,1) + \frac{1}{5}(1,1) + A_8$
        \item $A_1 + \frac{1}{3}(1,1) + \frac{1}{5}(1,2) + \frac{1}{22}(1,7)$
        \item $A_1 + \frac{1}{3}(1,1) + \frac{1}{5}(1,2) + \frac{1}{33}(1,13)$
        \item $A_1 + \frac{1}{3}(1,1) + \frac{1}{5}(1,2) + \frac{1}{43}(1,19)$
    \end{enumerate}
    \item $q \geq 20$ except the case $A_1 + \frac{1}{3}(1,1) + \frac{1}{5}(1,1)+ A_8$.
    \item  Let $\mathcal{D}$ be the reduced $f$-exceptional divisor on $S'$, and $E$ be any $(-1)$-curve on $S'$. Then,  $E.\mathcal{D} \geq 2$ and the equality holds if and only if  $E.f^{-1}(p_i) = 0$ for $i = 1, 2, 3$ and $E.f^{-1}(p_4) = 2$.
\end{enumerate} 
\end{lemma}

For later use, we give a formula for computing   $K^2_S$.
 
\begin{lemma}\cite[Section 3]{HK2013}\label{K^2}
Let $S$ be as in Notation \ref{notation}.  Then, we have
$$K^2_S = 9-L + \underset{p \in Sing(S)}{\sum} \Big(tr_p - 2l_p - 2 + \frac{q_{1, p}+q_{l, p}+2}{q_p}\Big).$$
\end{lemma}

\begin{proof}
Every surface quotient singular point is log terminal. So, up to numerical equivalence, we can write 
$$K_{S'} \underset{num}{\equiv} f^{*}K_S -
 \sum_{p \in Sing(S)}{\mathcal{D}_p}$$
where $\mathcal{D}_p = \sum(a_jA_j)$ is an effective
 $\mathbb{Q}$-divisor with $0 \leq a_j < 1$  supported on $f^{-1}(p)=\cup A_j$   for each singular point $p$.  Then, 
$$K^2_S = K^2_{S'} +\sum_{p}{\mathcal{D}_pK_{S'}}  = K^2_{S'} - \sum_{p}{\mathcal{D}_p^2}.$$
Since  $K^2_{S'} = 12 - (3+L) = 9 - L$ by Noether's formula and for each $p$
$$\mathcal{D}^2_p = 2l_p - tr_p  + 2 - \dfrac{q_{1,p} +q_{l,p} + 2}{q_p}$$ 
by \cite[Lemma 3.1 (3)]{HK2013}, the result follows.
\end{proof}

\begin{corollary}\label{235K^2}
Let $S$ be a rational $\mathbb{Q}$-homology projective plane  of log general type with exactly four cyclic quotient singular points $p_1, p_2, p_3, p_4$ of orders $(2, 3, 5, q)$. Assume that $p_2$ is of type $\frac{1}{3}(1,1)$. Then, we have
\begin{enumerate}
    \item If $p_3$ is of type $A_4$, then either $tr = 3l-2$ and $K^2_S = -\frac{2}{3}+\frac{q_1+q_l+2}{q}$ or $tr = 3l-3$ and $K^2_S = -\frac{5}{3}+\frac{q_1+q_l+2}{q}$. 
    \item If $p_3$ is of type $\frac{1}{5}(1,2)$, then either $tr = 3l-4$ and $K^2_S = -\frac{4}{15}+\frac{q_1+q_l+2}{q}$ or $tr = 3l-5$ and $K^2_S = -\frac{19}{15}+\frac{q_1+q_l+2}{q}$. 
    \item If $p_3$ is of type $\frac{1}{5}(1,1)$, then either $tr = 3l-7$ and $K^2_S = -\frac{13}{15}+\frac{q_1+q_l+2}{q}$ or $tr = 3l-8$ and $K^2_S = -\frac{28}{15}+\frac{q_1+q_l+2}{q}$. 
\end{enumerate}
In the above, the subscript $\mathcal{D}_4$ is omitted for $tr$, $l$, $q_1$, $q_l$ and $q$. 
\end{corollary}
 
Recall that the orbifold Euler number is defined as follows: 
$$e_{orb}(S) = e_{top}(S) - \underset{p \in Sing(S)}{\sum} \Big (1-\frac{1}{|G_p|}\Big)$$
where $|G_p|$ denotes the order of the local fundamental group of $p$. Then, we have  

\begin{theorem}  \cite{Megyesi}   \label{bmy-original} 
Let $S$ be a normal projective surface with quotient singularities. If $K_S$ is nef, then  $$  K^2_S \leq 3e_{orb}(S).$$
\end{theorem}

\begin{corollary}\label{bmy}
Let $S$ be as in Lemma \ref{235type}. Then, 
$0 < K^2_S \leq \frac{3}{q}+\frac{1}{10}. $ 
In particular, $0 < K^2_S  \leq \frac{1}{4}.$ 
\end{corollary}

\begin{proof}
Since $K_S$ is ample, $K^2_S > 0$. Now the result follows from Lemma \ref{235type} (2) and Lemma \ref{235K^2} (3).
\end{proof}

As an application, we can prove a nonexistence of a $\mathbb{Q}$-homology projective plane with special  singular points for later use.

\begin{corollary} 
There exists no rational $\mathbb{Q}$-homology projective plane of log general type with   four cyclic quotient singular points of type $A_1 + \frac{1}{3}(1,1)+\frac{1}{5}(1,1)+ \frac{1}{2l+1}(1,l)$ for any $l \geq 1$.  
\end{corollary}

\begin{proof}
By Lemma \ref{235type} (1), $L \geq 11$, so $l \geq 8$. 
Since the singularity $\frac{1}{2l+1}(1,l)$ has the resolution graph of the form $\overset{-3}{\circ}-\overset{-2}{\circ}-\overset{-2}{\circ}-\ldots-\overset{-2}{\circ}$,  by Lemma  \ref{K^2}, $K^2_S <0$, a contradiction to  Lemma \ref{bmy}.  
\end{proof}

\begin{corollary}\label{L=13}
Let $S$ be a rational $\mathbb{Q}$-homology projective plane with   four cyclic quotient singular points of type $A_1 + \frac{1}{3}(1,1)+A_4+  [2,a,2,b,2,c,2]$ where the last one is the singularity type for the singular point whose dual graph is of the form $\overset{-2}{\circ}-\overset{-a}{\circ}-\overset{-2}{\circ}-\overset{-b}{\circ}-\overset{-2}{\circ}-\overset{-c}{\circ}-\overset{-2}{\circ}$ where $a,b$ and $c$ are integers greater than $1$. Assume that $a+b+c=10$ or  $a+b+c=11$.  Then, $S$ is not of log general type.
\end{corollary}

\begin{proof}
Following Notation \ref{notation}, we have 
$q_1 = 8abc-8ab-4bc-8ca+6a+4b+2c-1,$ $q_l = 8abc-4ab-8bc-8ca+2a+4b+6c-1$ and $q= 16abc-16ab-16bc-16ca+12a+16b+12c-8$. Thus, we see that
$$\frac{q_1+q_l+2}{q} -1 = \frac{4(b-1)(a+c-2)}{q}.$$
If $a+b+c=11$, then by Lemma \ref{K^2}, Corollary \ref{235K^2} and Corollary \ref{bmy}, we have 
$$\frac{8}{q} + \frac{1}{3} \leq K^2_S = \frac{4(b-1)(a+c-2)}{q} + \frac{1}{3} \leq \frac{3}{q}+\frac{3}{20}, $$
which is a contradiction. If $a+b+c=10$, then by Lemma \ref{K^2} and Corollary \ref{235K^2}, we have 
$$ K^2_S = \frac{4(b-1)(a+c-2)}{q} - \frac{2}{3}.$$
In this case, one can compute that $K^2_S <0$, which is a contradiction.
\end{proof}

\subsection{A curve-detecting formula}

We present a useful formula for detecting the existence of a $(-1)$-curve. 

\begin{proposition}\label{formular-original}\cite[Proposition 4.2]{HK2013}
 Let $C$ be a smooth irreducible curve   on $S'$.  Then, there exists a rational number $m_E$   satisfying both of the following.
    \begin{enumerate}
     \item $C.K_{S'} = m_CK^2_S - \underset{p}{\sum} \overset{ l_p}{\underset{j = 1}{\sum}} \big( 1 - \dfrac{v_{j,p} + u_{j,p}}{q_p} \big) CA_{j,p}.$
         \item 
If, for each $p \in                Sing(S)$, $C$ has a non-zero intersection number with at most $2$ components of $f^{-1}(p)$, i.e., $CA_{j,p}=0$ for $j\neq s_p, t_p$ for some $s_p$ and $t_p$ with $1\le s_p< t_p\le l_p$ where the rational curves $A_{1,p}, A_{2,p}, \ldots, A_{l,p}$ form the tree of exceptional divisor of $f^{-1}(p)$,
then\\
$$ C^2  =   m^2_C K^2_S - \underset{p}{\sum} \Big(\dfrac{v_{s_p} u_{s_p}}{q_p}(CA_{s_p})^2+ \dfrac{v_{t_p} u_{t_p}}{q_p} (CA_{t_p})^2                + \dfrac{2 v_{t_p} u_{s_p}}{q_p} (CA_{s_p})(CA_{t_p}) \Big) .$$
    \end{enumerate}
\end{proposition}

\begin{remark}
We omit the subscript $C$ of $m_C$ if there is no confusion in the context.
\end{remark}

One can determine the positivity of $K_S$ simply by looking at the sign of $m$ for an irreducible  curve $C$ on $S'$.

\begin{lemma}\cite[Lemma 2.6]{HK2012} \label{ample}
      \begin{enumerate}
      \item $K_S$ is ample 
      iff $m_C>0$ for all irreducible curves $C$ not contracted by $f$
      iff $m_C>0$ for an irreducible curve $C$ not contracted by $f$.
      \item $K_S$ is numerically trivial 
      iff $m_C=0$ for all irreducible curves $C$ not contracted by $f$
      iff $m_C=0$ for an irreducible curve $C$ not contracted by $f$. 
      \item $-K_S$ is ample 
      iff $m_C<0$ for all irreducible curves $C$ not contracted by $f$ 
      iff $m_C<0$ for an irreducible curve $C$ not contracted by $f$. 
      \end{enumerate}
 
\end{lemma}

We present some   applications of the formula that will be used later.

\begin{lemma}\label{23719final}
Let $S$ be a $\mathbb{Q}$-homology projective plane  with four singular points of type  $A_1+A_2+\frac{1}{7}(1,1)+\frac{1}{19}(1,9)$. Let $\mathcal{D}$ be the reduced $f$-exceptional divisor on $S'$, and $E$ be any $(-1)$-curve on $S'$.
Assume that $K_S$ is ample. 
Then,  $E.\mathcal{D} \geq 3$. 
\end{lemma}

\begin{proof}
By Lemma \ref{23719conf}, we may assume that $E.\mathcal{D} = E.\mathcal{D}_4 = 2$. Since $K_S$ is ample, by Lemma \ref{ample}, $m > 0$. By Proposition \ref{formular-original},  we have $$1 < \sum (1-\frac{v_j+u_j}{q}).$$ 
Since $E.\mathcal{D}_4 = 2$,  
this is impossible by the computation in Table \ref{23719}.

\begin{table}[ht]
\caption{}\label{23719}
\renewcommand\arraystretch{1.5}
\noindent\[
\begin{array}{|c|c|c|c|c|c|c|c|c|c|c|}
\hline
 & \multicolumn{9}{|c|}{[3,2,2,2,2,2,2,2,2]}  \\   \hline
j& 1 &2&3&4&5&6&7&8&9\\
 \hline
1-\frac{v_j + u_{j}}{q}   & \frac{9}{19}&\frac{8}{19}&  \frac{7}{19}&  \frac{6}{19}&  \frac{5}{19}&\frac{4}{19}&  \frac{3}{19}&  \frac{2}{19}&  \frac{1}{19}\\
\hline
\end{array}
\]
\end{table}

\end{proof}

\begin{lemma}\label{end-end}
Let $S$ be a $\mathbb{Q}$-homology projective plane with cyclic quotient singularities, and $S'$ be its minimal resolution. Assume that there exists a $(-1)$-curve $E$ on $S'$ such that the dual graph of $\mathcal{D}+E$ forms a cycle where $\mathcal{D}$ is the $f$-exceptional divisor of a singular point of $S$. Then, $K_S$ is not ample.
\end{lemma}

\begin{proof}
By Proposition \ref{formular-original}, we have 
$$m_EK^2 = 1 - \frac{q_1+q_l+2}{q} \qquad \text{ and  }  \qquad
0 \leq m^2_EK^2_S = \frac{q_1+q_l+2}{q} - 1.$$
Now Lemma \ref{ample} implies that $K_S$ is not ample. 
\end{proof}

\begin{lemma}\label{RDP-case}
Let $S$ be a $\mathbb{Q}$-homology projective plane  with exactly four cyclic quotient singular points $p_1, p_2, p_3, p_4$ of orders $(2, 3, 5, q)$ with $q \geq 2$.
\begin{enumerate}
    \item Assume that $K_S$ is ample. If $p_4$ is a rational double point, then the singularities of $S$ are of type $A_1 + \frac{1}{3}(1,1)+\frac{1}{5}(1,1)+A_8$ and    $E.\mathcal{D} \geq 3$ for any $(-1)$-curve $E$.     
    \item There exists no $(-1)$-curve $E$ with $E.\mathcal{D}=E.D_1=2$ where $D_1$ is an end component of $\mathcal{D}_4$.
\end{enumerate} 
\end{lemma}

\begin{proof}
(1)  If $p_3$ is of type $\frac{1}{5}(1,1)$, then $tr_{\mathcal{D}_4} \geq 3l_{\mathcal{D}_4}-8$ by Lemma \ref{235K^2}. But, since $p_4$ is a rational double point, $tr_{\mathcal{D}_4} = 2l_{\mathcal{D}_4}$. Thus, we have $l \leq 8$, hence $L \leq 11$. Then, by Lemma \ref{235type} (1), $p_4$ is of type $A_1 + \frac{1}{3}(1,1)+\frac{1}{5}(1,1)+A_8$. 
If $E.\mathcal{D} \leq 2$, then by Lemma \ref{235type},  $E.\mathcal{D}_4 = 2$. Since   $\mathcal{D}_4$  is an exceptional divisor of a  rational double point,  $m_E K^2_S = -1$ by Proposition \ref{formular-original}. By Lemma \ref{ample}, $-K_S$ is ample,  a contradiction.  If $p_3$ is not of type $\frac{1}{5}(1,1)$, then by a similar argument as before using Lemma \ref{235K^2}, we have $L = 8$ or $9$, a contradiction. 

(2) Assume that such a curve $E$ exists. Let $n_1 := -D^2_1$. By Proposition \ref{formular-original},  
$$ 0 < m_EK^2_S = \frac{4q_1}{q}-1 = \frac{4q_1-q}{q} = \frac{(4-n_1)q_1+q_{1,2}}{q},$$
hence $n_1 \leq 4$. By contracting $E$, we get a minimal resolution $\bar{S}'$ of another rational $\mathbb{Q}$-homology projective plane $\bar{S}$. We  use  the  same  notation  for  the new surface and the curves lying on it    with the bar above the letter. We do not put the bar for $q_1$, $q_l$ and $q$ for readability since it does not cause any confusion. Let $\bar{C} := \bar{D}_1$.

Consider the case $n_1=4$. Then, $\bar{C}^2 = 0$. By Proposition \ref{formular-original}, we have
$$  m_{\bar{C}} K^2_{\bar{S}} = 1-\frac{q_1+1}{q} = \frac{q-q_1-1}{q}
\quad \text{     and     } \quad
 m^2_{\bar{C}} K^2_{\bar{S}} = \frac{q_1}{q}.$$ 

We claim that $q-q_1-1 > 0$. If not, then  $q=q_1+1$, so $p_4$ is a rational double point. By Lemma \ref{235K^2}, $K^2_{\bar{S}} \neq 0$, so $K_S$ is not numerically trivial. But   since 
$m_{\bar{C}} K^2_{\bar{S}} = 0$, Lemma \ref{ample} gives a contradiction.

Now, we see that $$m_{\bar{C}} = \frac{q_1}{q-q_1-1} > 0,$$ thus 
 $K_{\bar{S}}$ is ample by Lemma \ref{ample}.   
 Note that since $q = 4q_1 - q_{1,2} > 3q_1$, so $q \geq 3q_1+1$, thus $q -q_1-1 \geq 2q_1$. Moreover, since $4q_1 > q$, $$K^2_{S'} = \frac{(q-q_1-1)^2}{qq_1} \geq \frac{4q_1^2}{qq_1} > 1,$$ a contradiction by Lemma \ref{bmy}.

Consider the case $n_1 = 3$. Then, $\bar{C}^2 = 1$.  By Proposition \ref{formular-original}, we have
$$m_{\bar{C}}K^2_{\bar{S}} = \Big(1-\frac{q_1+1}{q}\Big)-1 = -\frac{q_1+1}{q} \qquad \text{  and  } \qquad m^2_{\bar{C}}K^2_{\bar{S}} = 1+\frac{q_1}{q}=\frac{q+q_1}{q}.$$
So, $m = -\frac{q+q_1}{q_1+1}$, hence, by Lemma \ref{ample}, $\bar{S}$ is a log del Pezzo surface of Picard number one. Moreover, $$K^2_{\bar{S}} = \frac{(q_1+1)^2}{q(q_1+q)}.$$
By \cite[Lemma 3]{HK2011b}, the integer $30qK^2_{\bar{S}}$ is a square number, thus so is $30(q_1+q)$. Then, $gcd(q, 30) \neq 1$. Indeed, if otherwise, such a surface does not exist by  \cite{HK2014}, as a solution to the algebraic Montgomery-Yang problem for log del Pezzo surfaces of Picard number one. Thus, $q$ is divisible by $2$, $3$, or $5$. It is easy to see that none of them is possible since $gcd(q, q_1) = 1$. For instance, if $q$ is even, so is $q_1$ since $30(q_1+q)$ is a square. But this is a contradiction since $gcd(q,q_1) = 2 \neq 1$. 

Consider the case $n_1=2$. By Corollary \ref{235K^2}, it is not hard to see that $$\frac{q_1+q_l+2}{q} = \frac{\bar{q}_1+\bar{q}_l+2}{\bar{q}}. $$
We claim that this is equivalent for $p_4$ to be a rational double point. Then, since $L \geq 11$ by Lemma \ref{235type}, this leads to a contradiction by Corollary \ref{235K^2}. For example, if $p_3$ is of type $\frac{1}{5}(1,1)$, then since $tr = 3l-7$, we have $l =7$, hence $L \leq 10$, a contradiction by Corollary \ref{235type}.

Now we prove the claim. It is easy to see that $q_1 = \bar{q}$, $q_l = 2\bar{q}_{\bar{l}} -\bar{q}_{1,\bar{l}}$ and $q = 2\bar{q}-\bar{q}_1$. Thus, 
$$\frac{q_1+q_l+2}{q} - \frac{\bar{q}_1+\bar{q}_{\bar{l}}+2}{\bar{q}}=\frac{\bar{q}+2\bar{q}_{\bar{l}}-\bar{q}_{1,\bar{l}}+2}{2\bar{q}-\bar{q}_1}-\frac{\bar{q}_1+\bar{q}_{\bar{l}}+2}{\bar{q}}.$$
Multiplying $\bar{q}(2\bar{q}-\bar{q}_1)$, by using the well-known identity  $\bar{q}_1\bar{q}_{\bar{l}} = \bar{q}\bar{q}_{1,\bar{l}}+1$, we can easily see that the above expression becomes
$$(\bar{q}-\bar{q}_1-1)^2$$
which equals to $0$ if and only if $p_4$ is a rational double point.
\end{proof}

\section{Singular fibers of a $\mathbb{P}^1$-fibration} 

We study singular fibers of a  $\mathbb{P}^1$-fibration on the minimal resolution of some rational $\mathbb{Q}$-homology projective planes. 
See \cite{Miyanishi} or \cite{GMM} for basics about $\mathbb{P}^1$-fibrations on rational surfaces.

In this section, we always denote by $S'$  the minimal resolution of a rational $\mathbb{Q}$-homology projective plane of log general type with $4$ cyclic singularities of orders $(2,3,5,q)$ where $q \geq 20$  and the order $3$ singularity is of type $\frac{1}{3}(1,1)$.  Assume that $S'$ admits a $\mathbb{P}^1$-fibration $\Phi: S' \rightarrow \mathbb{P}^1$ that has at most four horizontal components; at most one of them being a $2$-section and the rests being ordinary sections. Assume further that $\mathcal{D}_1$, $\mathcal{D}_2$ and $\mathcal{D}_3$ belong to fiber components of $\Phi$.

\begin{lemma}\label{P^1-fib}
Let $F$ be a singular fiber of $\Phi$ containing a connected component of $\mathcal{D}$, say $\mathcal{D}_k$, where $k \neq 4$. Let $E$ be a $(-1)$-curve in $F$ intersecting $\mathcal{D}_k$. Then, we have the following.
\begin{enumerate}
    \item $E.(F-E) \leq 2$.
    \item $E$ intersects a horizontal component of $\Phi$.
\end{enumerate}
\end{lemma}

\begin{proof}
(1) If $E.(F-E) \geq 3$, then $F^2 >0,$ a contradiction. 

(2) By Lemma \ref{235type}, $E.\mathcal{D} \geq 3$. By (1),   there  exists a component of $\mathcal{D}$  that is a horizontal component of $\Phi$ intersecting $E$. 

\end{proof}

\begin{lemma}\label{fiberD_2}
Let $F$ be the singular fiber of $\Phi$ containing the $(-3)$-curve $\mathcal{D}_2$. Then, $\Phi$ has a $2$-section and we have the following two cases for the configuration of $F$.
    \begin{enumerate}
        \item $F=E_1+A+2E_2+B$ whose dual graph is of the form   $\overset{-1}{\underset{E_1}{\circ}}-\overset{-3}{\underset{A}{\circ}}-\overset{-1}{\underset{E_2}{\circ}}-\overset{-2}{\underset{B}{\circ}}$
        where  $B$ is a component of $\mathcal{D}_1$ or $\mathcal{D}_4$.  
        \item $F=E_1+A+2E_2+B$ whose dual graph is of the form   $\overset{-1}{\underset{E_1}{\circ}}-\overset{-3}{\underset{A}{\circ}}-\overset{-1}{\underset{E_2}{\circ}}-B$
        with $B^2= -2$ where either $B = C_1 + 2E_3+C_2$ whose dual graph is of the form $\overset{-4}{\underset{C_1}{\circ}}-\overset{-1}{\underset{E_3}{\circ}}-\overset{-2}{\underset{C_2}{\circ}}$ where $C_1$ intersects $E_2$, or $B$ is obtained  by a finite sequence of blowups at the intersection point of two irreducible curves in the previous chain of $B$ one of them being a $(-1)$-curve.
    \end{enumerate}
In both cases, we have  $s_1.E_1=s_2.E_1=1$ and $s.E_2=1$ where $s_1$ and $s_2$ are sections of $\Phi$; and $s$ is a $2$-section of $\Phi$.
\end{lemma}

\begin{proof}
Let $t$ be the number of $(-1)$-curves in $F$ intersecting $\mathcal{D}_2$. 
Note that $1 \leq t \leq 3$ since otherwise we have $F^2 \neq 0$. If $t=1$, then $E$ has multiplicity $3$, thus it should intersect a horizontal component of $\Phi$ with multiplicity at least $3$, a contradiction.

If $t=3$, then $F$ is of the form $\mathcal{D}_2+E_1+E_2+E_3$ where $E_1$, $E_2$ and  $E_3$ are disjoint $(-1)$-curves intersecting the $(-3)$-curve $\mathcal{D}_2$. By Lemma \ref{235type} (3), each of $E_1$, $E_2$ and $E_3$ should intersect at least two more components of $\mathcal{D}$ that  are horizontal components of $\phi$. This is a contradiction since  $\Phi$ has at most $4$ horizontal components and at most one $2$-section.

Thus, $t = 2$. Now, the configuration of $F$ is of the form (2) in the statement. Note that  $E_1$ has multiplicity one and $E_2$ has multiplicity two. Moreover, $B^2=-2$.   If $B$ consists of one irreducible component, we arrive at  (1). In this case,  $B$ is not a component of $\mathcal{D}_3$ since every irreducible component of $\mathcal{D}_3$ belongs to one fiber. 
\end{proof}

From now on, assume that $S'$ admits a $\mathbb{P}^1$-fibration $\Phi: S' \rightarrow \mathbb{P}^1$ that has three horizontal components consisting of  two sections $s_1$ and $s_2$; and a  $2$-section $s$.

\begin{lemma}\label{P^1-fibration}
Let $F$ be a singular fiber of $\Phi$ containing a connected component $\mathcal{D}_k$ of $\mathcal{D}$ where $k \neq 4$.  
Let $E$ be a $(-1)$-curve in $F$ intersecting $\mathcal{D}_k$, and $t$ be the number of all such $(-1)$-curves. Then, we have the following.
\begin{enumerate} 
    \item We have $1 \leq t \leq 4$, and $t=4$ if and only if every $(-1)$-curve in $F$ has multiplicity one.
    \item The multiplicity of $E$ in $F$ is at most $2$.
    \item There exists at most one  $(-1)$-curve of multiplicity $2$ in $F$.
\end{enumerate}
\end{lemma}
   
\begin{proof}
(1) Since $F$ is a singular fiber of $\Phi$, $t \geq 1$. Since $\Phi$ has three horizontal components with one $2$-section, we have the result by Lemma \ref{P^1-fib} (2).

(2) It follows by Lemma \ref{P^1-fib} (2) since any horizontal component of $\Phi$ has multiplicity at most two.

(3) It follows by Lemma \ref{P^1-fib} (2) since $\Phi$ has exactly one $2$-section.
\end{proof}

Now we describe the configuration of possible singular fibers of $\Phi$ together with the horizontal components.

\begin{lemma}\label{fiberD_3}
Let $F$ be a singular fiber of $\Phi$ containing $\mathcal{D}_3$. Then, we have the following two cases.
    \begin{enumerate}
        \item $p_3$ is of type $A_4$ and $F=E_1+A_1+A_2+A_3+A_4+E_2$ whose dual graph is of the form   $\overset{-1}{\underset{E_1}{\circ}}-\overset{-2}{\underset{A_1}{\circ}}-\overset{-2}{\underset{A_2}{\circ}}-\overset{-2}{\underset{A_3}{\circ}}-\overset{-2}{\underset{A_4}{\circ}}-\overset{-1}{\underset{E_2}{\circ}}$. Here, $s_1.E_1=s_2.E_1=1$ and $s.E_2=2$. 
        \item $p_3$ is of type $\frac{1}{5}(1,2)$ and $F=E_1+A_1+A_2+2E_2+B$ whose dual graph is of the form   $\overset{-1}{\underset{E_1}{\circ}}-\overset{-2}{\underset{A_1}{\circ}}-\overset{-3}{\underset{A_2}{\circ}}-\overset{-1}{\underset{E_2}{\circ}}-\overset{-2}{\underset{\mathcal{D}_1}{\circ}}$. Here, $s_1.E_1=s_2.E_1=1$ and $s.E_2=1$. 
    \end{enumerate} 
\end{lemma}

\begin{proof}
Let $t$ be the number of $(-1)$-curves in $F$ intersecting $\mathcal{D}_3$.

Assume that $p_3$ is of type $\frac{1}{5}(1,1)$. By Lemma \ref{P^1-fibration} (1), $t \leq 4$. If $t=1$, then the unique $(-1)$-curve intersecting $\mathcal{D}_3$ has multiplicity $5$, a contradiction to Lemma \ref{P^1-fibration} (2). Similar argument shows that $t \neq 2$. If $t=3$, then the three $(-1)$-curves intersecting $\mathcal{D}_3$ have multiplicity $1$, $2$ and $2$, respectively. This gives a contradiction by Lemma \ref{P^1-fibration} (3). If $t=4$, then similar argument leads to a contradiction by  Lemma \ref{P^1-fibration} (1). 

Assume that $p_3$ is of type $A_4$. Let $E$ be a $(-1)$-curve intersecting $$\mathcal{D}_3 = \overset{-2}{\underset{A_1}{\circ}}-\overset{-2}{\underset{A_2}{\circ}}-\overset{-2}{\underset{A_3}{\circ}}-\overset{-2}{\underset{A_4}{\circ}}.$$ Then,  $E.A_2 = E.A_3 = 0$ since if otherwise we have $F^2 > 0$. Similarly,  we have $t \leq 2$. If $t=1$, then $E$ has multiplicity $5$, a contradiction to Lemma \ref{P^1-fibration} (2). If $t=2$, then we arrive at the case (1) in the statement.

Assume that $p_3$ is of type $\frac{1}{5}(1,2)$. Then, $\mathcal{D}_3 = \overset{-2}{\underset{A_1}{\circ}}-\overset{-3}{\underset{A_2}{\circ}}$. By Lemma \ref{P^1-fibration} (1), $t \leq 4$. Furthermore, $t \neq 4$ since $F^2 = 0$.  If $t=1$, then the unique $(-1)$-curve in $F$ has multiplicity $5$, a contradiction to Lemma \ref{P^1-fibration} (2).

Assume that $t=3$. Let $E_1$, $E_2$ and $E_3$ be the $(-1)$-curves in $F$. Then, the support of $F$ consists of $E_1, E_2, E_3$ and $\mathcal{D}_3$. Moreover,  $E_i$ has multiplicity one for each $i=1,2,3$.  This is a contradiction by Lemma \ref{P^1-fib} since $\Phi$ has $3$ horizontal components with only one $2$-section.

Assume that $t=2$. Then, $F$ is of the form $\overset{-1}{\underset{E_1}{\circ}}-\overset{-2}{\underset{A_1}{\circ}}-\overset{-3}{\underset{A_2}{\circ}}-\overset{-1}{\underset{E_2}{\circ}}-\overset{}{\underset{B}{\Box}}$ where $B$ is a tree of rational curves with $B^2=-2$. By Lemma \ref{235type} (3), we see that $E_1$ intersects two sections $s_1$ and $s_2$, and $E_2$ intersects the $2$-section $s$. Thus, $B$ does not contain a component of $\mathcal{D}_4$. Hence, $B=\mathcal{D}_1$. 
\end{proof}

\section{Proof of Theorem \ref{maintheorem}}

Under the situation in Notation \ref{notation},    we further assume that $S$ is of log general type and $H_1(S^{sm}, \mathbb{Z}) = 0$ where $S^{sm}$ denotes the smooth locus of $S$. Then,  we may assume that  the orders of the local fundamental groups of the singular points of $S$ is $(2,3,5,q)$ for some integer $q$ with $gcd(q,30)=1$ and the order $3$ singular point is of type $\frac{1}{3}(1,1)$  by Lemma  \ref{first4-tuple} and Lemma \ref{23719final}. By Lemma \ref{235type} and Lemma \ref{RDP-case}, we may assume that $L \geq 12$ and $q \geq 20$ except for the cases (b), (c) and  (d)  in Lemma \ref{235type}. Let $D_1,\ldots, D_l$ be irreducible components of $\mathcal{D}_4$ that form   a  chain of rational curves in this order.

Let $E$ be a $(-1)$-curve on $S'$ with $E.\mathcal{D} \leq 2$. By Corollary  \ref{ED>=2}, we may assume that   $E.\mathcal{D} = 2$. By Lemma \ref{235type}, $E$ intersects two components $D_i$ and $D_k$ of $\mathcal{D}_4$ with $i \leq k$.

\subsection{The case $i = k$}  The case $i=k=1$ is considered in Lemma \ref{RDP-case} (2) in a slightly more general context. By symmetry, we may assume that $1 < i=k < l$. Then, by contracting $E$ on $S'$, we get a minimal resolution of a rational $\mathbb{Q}$-homology projective plane with $5$ quotient singularities, a contradiction to Theorem  \ref{maximum}.

\subsection{The case $i < k$} By \cite[Lemma 4.1]{Zhang}, which can also be  proven using Proposition \ref{formular-original} as mentioned in \cite[Lemma 3.2]{HK2014},  we see that    $D^2_i = -2$ or $D^2_k = -2$. 
\subsubsection{The case $1 = i < k  \leq l$} By Lemma \ref{end-end}, we may assume that $k < l$. If $D^2_1 \leq -3$ and $D^2_k = -2$, then  by contracting $E$, we get a minimal resolution of a rational $\mathbb{Q}$-homology projective plane with $5$ quotient singularities, a contradiction to Theorem \ref{maximum}. This shows that $D^2_1 = -2$.

Assume that $D^2_k < -2$. Contracting $E$ and then the new $(-1)$-curves, if they exist, which are images of the components of $\mathcal{D}_4$, we arrive at the following three cases: (1) $D^2_1 = D^2_k = -2$, (2) $i=k=1$, or (3) $D^2_1 \leq -3$ and $D^2_k = -2$. Here, we need some care since the new $\mathbb{Q}$-homology projective plane might not be of log general type. Case (3) is impossible by the same argument as in the beginning of 4.2.1, and Case (2) was treated in Lemma \ref{RDP-case} (2).

Thus, it remains to consider the case  $D^2_1 = D^2_k = -2$. 
If $k =2$, then by contracting $E$ and the image of $D_1$, we get a minimal resolution of another rational $\mathbb{Q}$-homology projective plane $T$. Since $D^2_1=D^2_2= -2$, we see that $K^2_T <0$ by Lemma \ref{K^2} and Corollary \ref{235K^2}, a contradiction. 
If $k > 2$, then $D_1+2E+D_k$ induces a $\mathbb{P}^1$-fibration $\Phi: S' \rightarrow \mathbb{P}^1$ which has only ordinary sections. Consider the fiber $F$ containing the $(-2)$-curve $\mathcal{D}_2$. By Lemma \ref{fiberD_2}, $\Phi$ has a $2$-section, a contradiction.

\subsubsection{The case $1 < i < k  < l$} It is easy to see that $D^2_i = D^2_k = -2$ as in the beginning of 4.2.1. We claim that if $k = i+1$, then $p_4$ is a rational double point, i.e., it is of type $A_l$. This leads to  a contradiction by Lemma \ref{RDP-case}. To prove the claim, we first want to show that $D^2_{i-1} = -2$. If otherwise,   contracting $E$ and then the image of $D_i$, we get a minimal resolution of a rational $\mathbb{Q}$-homology projective plane with $5$ quotient singularities, a contradiction to Theorem \ref{maximum}.   By a similar argument, we   see   that $D^2_1= D^2_2 = \ldots = D^2_{i-2} = -2.$ By symmetry, we also have 
$D^2_{k+1}= D^2_{k+2} = \ldots = D^2_{l} = -2.$

Thus, we have $k \geq i+2$. Then, $D_i+2E+D_k$ induces a $\mathbb{P}^1$-fibration $\Phi: S' \rightarrow \mathbb{P}^1$. Since $D_i+2E+D_k$ forms one singular fiber of $\Phi$, we see that there are at most $4$ horizontal components, each of them is a section or possibly at most one double section of the $\mathbb{P}^1$-fibration $\Phi$.  Note that the other three connected components $\mathcal{D}_1 + \mathcal{D}_2 + \mathcal{D}_3$ form part of the fiber components. Take a fiber $F$ containing $\mathcal{D}_2$.  By Lemma \ref{fiberD_2},   $\Phi$ has a $2$-section. Then, the $2$-section is $D_{i+1}$ and $k=i+2$. Thus, $\Phi$ has   exactly three horizontal components, two of them being ordinary sections and the remaining one being a $2$-section. 

Since $p_3$ is not of type $\frac{1}{5}(1,1)$ by Lemma \ref{fiberD_3}, $F$ has a component of $\mathcal{D}_1$ or $\mathcal{D}_4$ by Lemma \ref{fiberD_2}. Since the $(-2)$-curve in $F$ does not intersect a horizontal component, $F$ does not contain a component of $\mathcal{D}_4$, hence $F$ is of the form in Lemma \ref{fiberD_2} (1) with $B = \mathcal{D}_1$. Let $G$ be the singular fiber of $\Phi$ containing $\mathcal{D}_3$. Then, since  $\mathcal{D}_1$ and $\mathcal{D}_2$ belong to the same singular fiber, we see that $F$ is of the form in Lemma \ref{fiberD_3} (1) by Lemma \ref{end-end}. In particular, $p_3$ is of type $A_4$.  By Lemma \ref{235type} (1),  $l \geq 6$. 
 
It remains to analyze the singular fibers with components from $\mathcal{D}_4$. Since there is only one multiple section that is a $2$-section,  $\mathcal{D}_4$ is of the form
$$\underset{D_1}{\overset{-2}{\circ}}-\underset{D_2}{\overset{-a}{\circ}}-\underset{D_3}{\overset{-2}{\circ}}-\underset{D_4}{\overset{-b}{\circ}}-\underset{D_5}{\overset{-2}{\circ}}-\underset{D_6}{\overset{-c}{\circ}}-\underset{D_7}{\overset{-2}{\circ}}$$
 for some integers $a,b,c \geq 2$ with $a+b+c=10$ or $11$. Note that $D_1+2E+D_7$ forms the remaining singular fiber of $\Phi$ for some $(-1)$-curve $E$, and we may regard $s_1=D_2$, $s_2=D_6$ and $s=D_4$.  This cannot happen by Corollary \ref{L=13}.  This completes the proof.

\bigskip {\bf Acknowledgements.}  
This research was supported by the Samsung Science and Technology Foundation under Project SSTF-BA1602-03.
 
\bigskip


\end{document}